\documentclass[12pt,reqno]{amsart}
\usepackage[english]{babel}
\usepackage{amsmath,amsfonts,amssymb,amsthm,amscd,latexsym}
\usepackage[T1]{fontenc}
\usepackage[utf8]{inputenc}

\usepackage{tikz}
\usetikzlibrary{arrows}
\usepackage{color}
\usepackage{cite}
\usepackage{multirow}
\usepackage{cite}

\newtheorem{theorem}{Theorem}[section]

\newtheorem{lemma}[theorem]{Lemma}

\newtheorem{definition}[theorem]{Definition}

\newtheorem{rem}[theorem]{Remark}

\numberwithin{equation}{section}
\DeclareFontFamily{OT1}{pzc}{}
\DeclareFontShape{OT1}{pzc}{m}{it}{<-> s * [1.10] pzcmi7t}{}
\DeclareMathAlphabet{\mathpzc}{OT1}{pzc}{m}{it}

\def\ms{\mathfrak{s}}

\def\mh{\mathfrak{h}}

\def\a{\alpha}

\def\C{\Bbb C}

\begin{document}

\title{Local derivation on the $n$-th Schr{\"o}dinger algebra}

%\begin{center}	
%{\bf Amir Alauadinov\footnote{V.I.Romanovskiy Institute of Mathematics Academy of Science of Uzbekistan; Ch.Abdirov 1, Department of Mathematics, Karakalpak State University, Nukus 230113, Uzbekistan; \ {amir\_t85@mail.ru}}},
%{\bf Doston Jumaniyozov\footnote{V.I.Romanovskiy Institute of Mathematics Academy of Science of Uzbekistan; \ {dostondzhuma@gmail.com}}},
%{\bf
%Bakhtiyor Yusupov\footnote{ V.I.Romanovskiy Institute of Mathematics Academy of Science of Uzbekistan; \
%Department of Physics and Mathematics, Urgench State University;
% \
%baxtiyor\_yusupov\_93@mail.ru}}
%\end{center}

\author[A.K.Alauadinov]{Amir Alauadinov$^{1,2}$}

\address{$^1$V.I.Romanovskiy Institute of Mathematics Uzbekistan Academy of Sciences. 9,  University street, 100174, Tashkent, Uzbekistan}

\address{$^2$Karakalpak State University, Ch.Abdirov 1, Nukus 230113, Uzbekistan}

\address{$^3$Department of Physics and Mathematics, Urgench State University, H. Alimdjan street, 14, Urgench
220100, Uzbekistan}

\author[B.B.Yusupov]{Bakhtiyor Yusupov$^{1,3}$}

\author[D.E.Jumaniyozov]{Doston Jumaniyozov$^1$}

\email{{\tt amir\_t85@mail.ru, baxtiyor\_yusupov\_93@mail.ru, dostondzhuma@gmail.com}}

%\address{  \newline \indent National University of Uzbekistan} \email{{\tt }}

\begin{abstract} This paper is devoted to study local derivations on the $n$-th Schr{\"o}dinger algebra $\mathcal{S}_{n}.$ We prove that every local derivation on $\mathcal{S}_{n}$  is a derivation.
\end{abstract}


\maketitle


\section{Introduction}\label{sec:intro}

Local derivations and automorphisms have been objects of many research works in recent years. Even many results have been addressed in this direction, there are still lots of related unsolved problems. R.V.Kadison, D.R.Larson and A.R.Sourour first introduced the notion of local derivations on algebras in their remarkable paper \cite{Kad, Lar}. After that many research papers have been published regarding local derivations of non-associative algebras (e.g. see \cite{AyuKud, AyuKhu, KudOmiKur, AlaYus1}). For instanxe, it is proved that every local derivation  on a finite-dimensional semisimple Lie algebra $\mathcal{L}$ over an algebraically closed field of characteristic zero is a derivation (see \cite{AyuKud}). Moreover, local derivations of solvable Lie algebras are investigated in \cite{AyuKhu}, and it is shown that  there exist solvable Lie algebras which admit local derivations which are not ordinary derivation. However, it is also proved that there are solvable Lie algebras for which every local derivation is a derivation (see \cite{AyuKhu}). Recently, in \cite{KudOmiKur}, the authors proved that every local derivation on solvable Lie algebras whose nilradical has maximal rank is a derivation. Further, it also is proved  that every local derivation on the conformal Galilei algebra is a derivation (see \cite{AlaYus1}). Recently, in \cite{AlaYus2,AlaYus3}, the authors proved that every local derivation on the Schrödinger algebra $\mathfrak{s}_n$ in $(n+1)$-dimensional space-time is a derivation when $n \in\mathbb{N}$. A.Alauadinov and B.Yusupov proved similar results
concerning local super-derivations on the super Schr\"{o}dinger algebras in their recent paper \cite{AlaYus4}. Note that Lie algebras that are close to nilpotent admits pure local derivations. Therefore, it is natural to conjecture that any nilpotent Lie algebra with nilindex greater than $3$ admits pure local derivations. In \cite{KJ}, it is proved that any $2$-step nilpotent Lie algebra over a field of arbitrary characteristic except $2,$ admits a local derivation which is not derivation.

We should note that all aforementioned results have been obtained in finite-dimensional algebras. In the infinite-dimensional case,
it is proved that every local derivation on some class of the locally simple Lie algebras, generalized Witt algebras, Witt algebras, the Witt algebras over a field of prime characteristic is a derivation (see \cite{Kaim,AyuKudYus,Yao}).

There are many applications of non-semisimple Lie algebras and groups in physics. One of the well-known non-semisimple Lie algebras is the Schrödinger algebra $\mathfrak{s}_n.$ The Schrödinger algebra is a Lie algebra of the Schrödinger group that is group of symmetries of the free particle Schrödinger equation (see \cite{Per,DDM1}). Recently there have been published a series of papers on studying the structure and representation theory of the Schrödinger algebra $\mathfrak{s}_1$ in the case of $(1+1)$-dimensional space-time, see \cite{AD, YT, ZCH}.	
	 	
In this paper, we study local derivations on $n$-th Schr{\"o}dinger algebra $\mathcal{S}_n$. Namely, we prove that every local derivation on $\mathcal{S}_n$ is a derivation.

\section{Preliminaries}

In this section, we first recall the definition of $\mathcal{S}_n$. We know that  the general linear Lie algebra $\mathfrak{gl}_{2n}$  has  the natural representation on $\C^{2n}$ by left matrix  multiplication. Let $\{u_k, v_k\}, \ 1\leq k\leq n$ be a basis of $\C^{2n}$.

The Heisenberg algebra $\mh_{n}=\C^{2n}\oplus \C z$ is a $2$-step nilpotent Lie algebra with Lie bracket given by
$$[u_k, v_{k}]=z, \qquad [z, \mh_n]=0,\quad 1\leq k\leq n.$$

Recently, the $n$-th Schrödinger algebra  $\mathcal{S}_{n}$ was introduced \cite{Ta}.

\begin{definition}
The $n$-th Schr\"{o}dinger algebra $\mathcal{S}_{n}$ is a Lie algebra with a $\mathbb{C}$-basis
$$
\{e,f,h,z,u_{k},v_{k},\mid k=1,...,n\}
$$
equipped with the following non-trivial commutation relations
\begin{align*}\label{Bracket}
[h,e]&=2e, & [h,f]&=-2f, & [e,f]&=h,\\
[u_k, v_j]&=\delta_{kj}z, & [h,u_k]&=u_k, & [h,v_k]&=-v_k,\\
[e,v_k]&=u_k, & [f,u_k]&=v_k, & [z,\mathcal{S}_{n}]&=0,
\end{align*}
where $\delta_{kj}$ is the Kronecker Delta defined as $1$ for $k=j$ and as $0$ otherwise.
\end{definition}

 From the definition one can notice that the $n$-th Schr{\"o}dinger algebra $\mathcal{S}_{n}$ is a finite-dimensional, non-semisimple and non-solvable Lie algebra, and it is the semidirect product $3$-dimensional simple Lie algebra and $(2n+1)$-dimensional Heisenberg algebra, that is
$$\mathcal{S}_{n}=(\mathfrak{sl}_{2}\ltimes\mathfrak{h}_{n})$$
where $\mathfrak{sl}_{2}={\rm Span}_{\mathbb{C}}\{e,f,h\}$ is the 3-dimensional simple Lie algebra and $\mathfrak{h}_{n}={\rm Span}_{\mathbb{C}}\{z,u_{k},v_{k}\mid1\leqslant k\leqslant n\}$ is the Heisenberg Lie algebra.

A derivation on a  Lie algebra $\mathcal{L}$ is a linear map  $D: \mathcal{L} \rightarrow \mathcal{L}$  which satisfies the Leibniz rule:
\begin{equation}\label{Defder}
D([x,y]) = [D(x), y] + [x, D(y)], \quad \text{for any} \quad x,y \in \mathcal{L}.
\end{equation}
The set of all derivations of $\mathcal{L}$ is denoted by $\mathrm{Der}(\mathcal{L})$ and with respect to the commutation operation is a Lie algebra.

For any element $y\in \mathcal{L}$ the
left multiplication operator $\mathrm{ad}_x: \mathcal{L} \to \mathcal{L}$, defined as $\mathrm{ad}_x(y)=[x,y]$ is a derivation, and derivations of this form are called inner derivations.
The set of all inner derivations of $\mathcal{L},$ denoted by $\mathrm{Inn}(\mathcal{L}),$ is an ideal in $\mathrm{Der}(\mathcal{L}).$
\begin{definition}
A linear operator $\Delta$ is called a local derivation if for any $x \in \mathcal{L},$ there exists a derivation $D_x: \mathcal{L} \rightarrow \mathcal{L}$ (depending on $x$) such that
$\Delta(x) = D_x(x).$ The set of all local derivations on $\mathcal{L}$ we denote by $\mathrm{LocDer}(\mathcal{L}).$
\end{definition}

Let $n\geq 2.$ For any $1\leq l<k\leq n,$ one can easily verify that the linear map $\sigma_{lk}:\mathcal{S}_n\rightarrow \mathcal{S}_{n}$
defined below is a derivation:
\begin{equation}
    \begin{split}
        \sigma_{lk}(h)&=\sigma_{lk}(e)=\sigma_{lk}(f)=\sigma_{lk}(z)=0,\\
        \sigma_{lk}(u_j)&=\delta_{lj}u_k-\delta_{lj}u_l,\quad \sigma_{lk}(v_k)=\delta_{lj}v_k-\delta_{kj}v_l.\\
    \end{split}
\end{equation}

Clearly, $\sigma_{lk} (1\leq l<k\leq n)$ are outer derivations. We define another outer derivation $\tau$ as follows
\begin{equation}{\label{tau}}
    \begin{split}
       \tau(h)&=\tau(e)=\tau(f)=0,\ \tau(u_l)=\frac12u_l,\ \tau(v_l)=\frac12v_l,\ \tau(z)=z,\  1\leq l\leq n.
    \end{split}
\end{equation}

The following theorem is proved in \cite{BH}.
\begin{theorem}{\label{thmder}}
    $\mathrm{Der}(\mathcal{S}_n)=\mathrm{Inn}(\mathcal{S}_n)  \oplus \bigoplus\limits_{1\leq i<j\leq n} \mathbb{C}\sigma_{ij}\bigoplus\mathbb{C}\tau.$
\end{theorem}

It should be noted that, if $n=1,$ then $\mathrm{Der}(\mathcal{S}_n)=\mathrm{Inn}(\mathcal{S}_n)  \bigoplus \mathbb{C}\tau,$ where $\tau$ is defined as \eqref{tau}.

%\begin{equation}\begin{split}
%D(e)&=a_{f}h+2a_he-\sum\limits_{k=1}^{n}b_{k}u_k,\\
%D(h)&=-2a_ee+2a_ff-\sum\limits_{k=1}^na_ku_k-\sum\limits_{k=1}^{n}b_{k}v_k,\\
%D(f)&=a_{e}h-2a_{h}f-\sum\limits_{k=1}^{n}a_{k}v_{k},\\
%D(z)&=cz,\\
%D(u_j)&=a_fv_j+\left(\a_h+\frac{1}2c\right)u_j-b_jz-\sum\limits_{ 1\leq k<j}d_{k,j}u_k+\sum\limits_{ j<l\leq n}d_{{j,l}}u_l\\
%D(v_j)&=a_eu_j+\left(\frac{1}2c-a_h\right)v_j-a_jz-\sum\limits_{ 1\leq k<j}d_{{k,j}}v_k+\sum\limits_{ j<l\leq n}d_{j,l}v_l,\\ \end{split}\end{equation}

\section{Main results}

In this section, we will prove that every local derivation on the $n$-th Schr{\"o}dinger algebra $\mathcal{S}_n$ is a derivation.

\begin{theorem}\label{thm21}
Every local derivation on the $n$-th Schr{\"o}dinger algebra $\mathcal{S}_n,\ n \geq2$ is a derivation.
\end{theorem}

For any $x\in\mathcal{S}_n$ there exist
$$\alpha=\alpha_{e}e+\alpha_{h}h+\alpha_{f}f+\sum\limits_{k=1}^{n}\alpha_{u_k}u_k+\sum\limits_{k=1}^{n}\alpha_{v_k}v_k+\alpha_{z}z,$$
and $\lambda,\ \mu_{i,j}\in\mathbb{C}$ such that, by using Theorem \ref{thmder} we can write

$$\Delta(x)=[\alpha_x,x]+\lambda_x\sigma(x)+\sum\limits_{1\leq k< j\leq n}\mu_x\sigma_{k,j}(x),$$
where $\alpha_e,\alpha_f,\alpha_h,\alpha_z,\alpha_{u_k},\alpha_{v_k},\lambda,\ \mu_{k,j}$ are complex numbers depending on $x.$

Then we obtain
\begin{equation}\label{asos}
\begin{aligned}
\Delta(e)&=\alpha_f^{(e)}h+2\alpha_h^{(e)}e-\sum\limits_{k=1}^n\alpha^{(e)}_{v_k}u_k,\\
\Delta(h)&=-2\alpha_e^{(h)}e+2\alpha_f^{(h)}f-\sum\limits_{k=1}^n\alpha^{(h)}_{u_k}u_k-\sum\limits_{k=1}^n\alpha^{(h)}_{v_k}v_k,\\
\Delta(f)&=\alpha_e^{(f)}h-2\alpha_{h}^{(f)}f-\sum\limits_{k=1}^n\alpha^{(f)}_{u_k}v_k,\\
\Delta(u_j)&=\alpha^{(u_j)}_{f}v_j+\left(\alpha^{(u_j)}_{h}+\frac{\lambda^{(u_j)}}{2}\right)u_j-\alpha^{(u_j)}_{v_j}z-\sum\limits_{ 1\leq k<j}\mu^{(u_j)}_{{k,j}}u_k+\sum\limits_{ j<l\leq n}\mu^{(u_j)}_{{j,l}}u_l,\\
\Delta(v_j)&=\left(\frac{\lambda^{(v_j)}}{2}-\alpha^{(v_j)}_{h}\right)v_j+\alpha^{(v_j)}_{e}u_j+\alpha^{(v_j)}_{u_j}z-\sum\limits_{ 1\leq k<j}\mu^{(v_j)}_{{k,j}}v_k+\sum\limits_{ j<l\leq n}\mu^{(v_j)}_{{j,l}}v_l,\\
\Delta(z)&=\lambda^{(z)}z.
\end{aligned}
\end{equation}

For the proof of this Theorem we need several Lemmata.

In the following lemmas, we assume that $\Delta$ be a local derivation on $\ms_n$ such that $\Delta(h+z)=0.$

\begin{lemma}\label{lem2}
 $\Delta(h)=\Delta(z)=0.$
\end{lemma}
\begin{proof}
We consider,
$$0=\Delta(h+z)=\Delta(h)+\Delta(z)=-2\alpha_e^{(h)}e+2\alpha_f^{(h)}f-\sum\limits_{k=1}^n\alpha^{(h)}_{u_k}u_k-\sum\limits_{k=1}^n\alpha^{(h)}_{v_k}v_k+\lambda^{(z)}z,$$
which implies $\alpha_e^{(h)}=\alpha_f^{(h)}=\alpha_{u_k}^{(h)}=\alpha_{v_k}^{(h)}=\lambda^{(z)}=0,$ where $k\in\{1,2,\cdots,n\}.$ Then $\Delta(h)=\Delta(z)=0.$
\end{proof}

\begin{lemma}\label{lem3}

 Coefficients $a_f^{(e)},\ a_e^{(f)}$ $a^{(u_j)}_f$ and $a^{(v_j)}_{e}, 1\leq j\leq n$ in the formula (\ref{asos}) are equal to zero.
 %   \begin{equation}
%\begin{aligned}
%\Delta(e)&=2\alpha_h^{(e)}e-\sum\limits_{k=1}^n\alpha^{(e)}_{v_k}u_k,\\
%\Delta(f)&=-2\alpha_{h}^{(f)}f-\sum\limits_{k=1}^n\alpha^{(f)}_{u_k}v_k,\\
%\Delta(u_j)&=\left(\alpha^{(u_j)}_{h}+\frac{\lambda^{(u_j)}}{2}\right)u_j-\alpha^{(u_j)}_{v_j}z+\sum\limits_{ 1\leq k<j}\mu^{(u_j)}_{{k,j}}u_k-\sum\limits_{ j<l\leq n}\mu^{(u_j)}_{{j,l}}u_l,\\
%\Delta(v_j)&=\left(\frac{\lambda^{(v_j)}}{2}-%\alpha^{(v_j)}_{h}\right)v_j+\alpha^{(v_j)}_{u_j}z+\sum\limits_{ 1\leq k<j}\mu^{(v_j)}_{{k,j}}v_k-\sum\limits_{ j<l\leq n}\mu^{(v_j)}_{{j,l}}v_l.\\
%\end{aligned}
%\end{equation}
\end{lemma}

\begin{proof}
From the equalities, take an element $x=h+e$, we consider, using the rate of local derivation we calculate the following
\begin{equation*}\begin{split}
\Delta(x)&=\Delta(h+e)=[b,h+e]+c\sigma(h+e)+\sum\limits_{1\leq k< l\leq n}d_{k,l}\sigma_{k,l}(h+e)=\\
             &=\big[a_{f}f+a_{h}h+a_{e}e+a_zz+
             \sum\limits_{k=1}^na_ku_k+
             \sum\limits_{k=1}^nb_kv_k,\ h+e\big]+\\
             &+c\sigma(h+e)+\sum\limits_{1\leq k< l\leq n}d_{k,l}\sigma_{k,l}(h+e)=a_{f}h+2a_ff+\ast e-\sum\limits_{k=1}^{n}\ast u_k-\sum\limits_{k=1}^{n}\ast v_k.\\
\end{split}
\end{equation*}

On the other hand, based on (\ref{asos}), we calculate the following equality:
\begin{equation*}\begin{split}
\Delta(x)&=\Delta(h+e)=\Delta(h)+\Delta(e)=\alpha_f^{(e)}h+2\alpha_h^{(e)}e-\sum\limits_{k=1}^n\alpha^{(e)}_{v_k}u_k.
\end{split}
\end{equation*}
    Comparing the coefficients at the basis elements $h$ and $f,$ we get $a_f=0,\ a_f=\alpha_f^{(e)}$ which implies
$$\alpha_f^{(e)}=0.$$

  Consider element $x=h+f.$ Then by there exists a derivation $D_x$ with $\Delta(x)=D_x(x).$ So, we have the following
\begin{equation*}\begin{split}
\Delta(x)&=\Delta(h+f)=[b,h+f]+c\sigma(h+f)+\sum\limits_{1\leq k< l\leq n}d_{k,l}\sigma_{k,l}(h+f)=\\
             &=\big[a_{f}f+a_{h}h+a_{e}e+a_zz+
             \sum\limits_{k=1}^na_ku_k+
             \sum\limits_{k=1}^nb_kv_k,\ h+f\big]+\\
             &+c\sigma(h+f)+\sum\limits_{1\leq k< l\leq n}d_{k,l}\sigma_{k,l}(h+f)=
             -2a_ee+a_{e}h+\ast f-\sum\limits_{k=1}^n\ast u_k-\sum\limits_{k=1}^{n}\ast v_k.\\
\end{split}
\end{equation*}

On the other hand, based on (\ref{asos}), we calculate the following equality:
\begin{equation*}\begin{split}
\Delta(x)&=\Delta(h+f)=\Delta(h)+\Delta(f)=\alpha_e^{(f)}h-2\alpha_{h}^{(f)}f-\sum\limits_{k=1}^n\alpha^{(f)}_{u_k}v_k.
\end{split}
\end{equation*}
    Comparing the coefficients at the basis elements $h$ and $e,$ we get $a_e=0,\ a_e=\alpha_e^{(f)}$ which implies
$$\alpha_e^{(f)}=0.$$

 From the equalities, take an element $x=e+u_j$, for any fixed $1\leq j\leq n,$ we consider, using the rate of local derivation we have the following chain of equations:
\begin{equation*}\begin{split}
\Delta(x)&=\Delta(e+u_j)=[b,e+u_j]+c\sigma(e+u_j)+\sum\limits_{1\leq k< l\leq n}d_{k,l}\sigma_{k,l}(e+u_j)=\\
             &=a_{f}h+2a_he-\sum\limits_{k=1}^{n}b_{k}u_k+
             a_fv_j+\left(\a_h+\frac{1}2c\right)u_j-\\
             &-b_jz-\sum\limits_{ 1\leq k<j}d_{k,j}u_k+\sum\limits_{ j<l\leq n}d_{{j,l}}u_l.
\end{split}
\end{equation*}

On the other hand, based on (\ref{asos}), we get the following equality:
\begin{equation*}\begin{split}
\Delta(x)&=\Delta(e+u_j)=\Delta(e)+\Delta(u_j)=2\alpha_h^{(e)}e-\sum\limits_{k=1}^n\alpha^{(e)}_{v_k}u_k+\\
&+\alpha^{(u_j)}_{f}v_j+\left(\alpha^{(u_j)}_{h}+\frac{\lambda^{(u_j)}}{2}\right)u_j-\alpha^{(u_j)}_{v_j}z-\sum\limits_{ 1\leq k<j}\mu^{(u_j)}_{{k,j}}u_k+\sum\limits_{ j<l\leq n}\mu^{(u_j)}_{{j,l}}u_l.\\
\end{split}
\end{equation*}
    Comparing the coefficients at the basis elements $h$ and $v_j,$ we get $a_f=0,\ a_f=\alpha_f^{(u_j)}$ which implies
$$\alpha_f^{(u_j)}=0.$$

 From the equalities, take an element $x=f+v_j$, for any fixed $1\leq j\leq n,$ we consider, using the rate of local derivation we calculate the following
\begin{equation*}\begin{split}
\Delta(x)&=\Delta(f+v_j)=[b,f+v_j]+c\sigma(f+v_j)+\sum\limits_{1\leq k< l\leq n}d_{k,l}\sigma_{k,l}(f+v_j)=\\
             &=a_{e}h-2a_{h}f-\sum\limits_{k=1}^{n}a_{k}v_{k}+
             a_eu_j+\left(\frac{1}2c-a_h\right)v_j-\\
             &-a_jz-\sum\limits_{ 1\leq k<j}d_{{k,j}}v_k+\sum\limits_{ j<l\leq n}d_{j,l}v_l.
\end{split}
\end{equation*}

On the other hand, based on (\ref{asos}), we calculate the following equality:
\begin{equation*}\begin{split}
\Delta(x)&=\Delta(f+v_j)=\Delta(f)+\Delta(v_j)=-2\alpha_{h}^{(f)}f-\sum\limits_{k=1}^n\alpha^{(f)}_{u_k}v_k+\\
&+\left(\frac{\lambda^{(v_j)}}{2}-\alpha^{(v_j)}_{h}\right)v_j+\alpha^{(v_j)}_{e}u_j+\alpha^{(v_j)}_{u_j}z-\sum\limits_{ 1\leq k<j}\mu^{(v_j)}_{{k,j}}v_k+\sum\limits_{ j<l\leq n}\mu^{(v_j)}_{{j,l}}v_l.\\
\end{split}
\end{equation*}
    Comparing the coefficients at the basis elements $h$ and $u_j,$ we get $a_e=0,\ a_e=\alpha_e^{(v_j)}$ which implies
$$\alpha_e^{(v_j)}=0.$$
\end{proof}

\begin{lemma}\label{lem4}
Coefficients $\alpha^{(u_j)}_{v_j}$ and $\alpha^{(v_j)}_{u_j}, 1\leq j\leq n$ in the formula (\ref{asos}) are equal to zero.
  %  \begin{equation}
%\begin{aligned}
%\Delta(u_j)&=\left(\alpha^{(u_j)}_{h}+\frac{\lambda^{(u_j)}}{2}\right)u_j+\sum\limits_{ 1\leq k<j}\mu^{(u_j)}_{{k,j}}u_k-\sum\limits_{ j<l\leq n}\mu^{(u_j)}_{{j,l}}u_l,\\
%\Delta(v_j)&=\left(\frac{\lambda^{(v_j)}}{2}-\alpha^{(v_j)}_{h}\right)v_j+\sum\limits_{ 1\leq k<j}\mu^{(v_j)}_{{k,j}}v_k-\sum\limits_{ j<l\leq n}\mu^{(v_j)}_{{j,l}}v_l.\\
%\end{aligned}
%\end{equation}
\end{lemma}

\begin{proof}
 From the equalities, take an element $x=h+u_j$, for any fixed $1\leq j\leq n,$ we consider, using the rate of local derivation we calculate the following
\begin{equation*}\begin{split}
\Delta(x)&=\Delta(h+u_j)=[b,h+u_j]+c\sigma(h+u_j)+\sum\limits_{1\leq k< l\leq n}d_{k,l}\sigma_{k,l}(h+u_j)=\\
             &=2a_ff-b_jz-\sum\limits_{k=1}^{n}b_{k}v_k+a_fv_j
             -\ast e+\sum\limits_{k=1}^n\ast u_k.\\
\end{split}
\end{equation*}

On the other hand, based on (\ref{asos}), we calculate the following equality:
\begin{equation*}\begin{split}
\Delta(x)=\Delta(h+u_j)=\Delta(h)+\Delta(u_j)&=-\alpha^{(u_j)}_{v_j}z+\sum\limits_{k=1}^n\ast u_k.\\
\end{split}
\end{equation*}
    Comparing the coefficients at the basis elements $f,\ z$ and $v_j,$ we get $a_f=a_j=0,\ a_j=\alpha_{v_j}^{(u_j)}$ which implies
$$\alpha_{v_j}^{(u_j)}=0.$$

 From the equalities, take an element $x=h+v_j$, for any fixed $1\leq j\leq n,$ we consider, using the rate of local derivation we calculate the following
\begin{equation*}\begin{split}
\Delta(x)&=\Delta(h+v_j)=[b,h+v_j]+c\sigma(h+v_j)+\sum\limits_{1\leq k< l\leq n}d_{k,l}\sigma_{k,l}(h+v_j)=\\
             &=-2a_ee-a_jz-\sum\limits_{k=1}^na_ku_k+a_eu_j+\ast f-\sum\limits_{k=1}^{n}\ast v_k.\\
\end{split}
\end{equation*}

On the other hand, based on (\ref{asos}), we calculate the following equality:
\begin{equation*}\begin{split}
\Delta(x)=\Delta(h+v_j)=\Delta(h)+\Delta(v_j)&=\alpha^{(v_j)}_{u_j}z+\sum\limits_{k=1}^n\ast v_k.\\
\end{split}
\end{equation*}
    Comparing the coefficients at the basis elements $e,\ z$ and $u_j,$ we get $a_e=a_j=0,\ a_j=\alpha_{u_j}^{(v_j)}$ which implies
$$\alpha_{u_j}^{(v_j)}=0.$$
\end{proof}

\begin{lemma}\label{Lem5}
 Coefficients $\alpha_{u_j}^{(f)}$ and $\alpha^{(e)}_{v_j}, 1\leq i\leq n$ in the formula (\ref{asos}) are equal to zero and $\frac{\lambda^{(v_j)}}{2}-\alpha_h^{(v_j)}=-\alpha_h^{(e)},$ $\alpha_h^{(u_j)}+\frac{\lambda^{(u_j)}}{2}=\alpha_h^{(e)}.$
\end{lemma}

\begin{proof}
From the equalities, take an element $x=e+f$,

\begin{equation}\begin{split}\label{ef11}
\Delta(x)&=\Delta(e+f)=[b,e+f]+c\sigma(e+f)+\sum\limits_{1\leq k< l\leq n}d_{k,l}\sigma_{k,l}(e+f)=\\
&=a_{f}h+2a_he-\sum\limits_{k=1}^{n}b_{k}u_k+a_{e}h-2a_{h}f-\sum\limits_{k=1}^{n}a_{k}v_{k}.\\
\end{split}
\end{equation}

On the other hand, based on (\ref{asos}), we calculate the following equality:
\begin{equation}\begin{split}\label{ef22}
\Delta(x)&=\Delta(e+f)=\Delta(e)+\Delta(f)=2\alpha_h^{(e)}e-\sum\limits_{k=1}^n\alpha^{(e)}_{v_k}u_k-2\alpha_{h}^{(f)}f-\sum\limits_{k=1}^n\alpha^{(f)}_{u_k}v_k.\\
\end{split}
\end{equation}
    Comparing the coefficients at the basis elements $e$ and $f$, (\ref{ef11}) and (\ref{ef22}) we abtain that
\begin{equation*}
    \alpha_h^{(e)}=\alpha_h^{(f)}.
\end{equation*}

From the equalities, take an element $x=f-\frac{1}{2}z+v_j$, for any fixed $1\leq j\leq n,$

\begin{equation}\begin{split}\label{f1}
\Delta(x)&=\Delta(f-\frac{1}{2}z+v_j)=[b,f-\frac{1}{2}z+v_j]+c\sigma(f-\frac{1}{2}z+v_j)+\\
&+\sum\limits_{1\leq k< l\leq n}d_{k,l}\sigma_{k,l}(f-\frac{1}{2}z+v_j)=
             a_{e}h-2a_{h}f-\sum\limits_{k=1}^{n}a_{k}v_{k}-\frac{c}{2}z+\\
             &+a_eu_j+\left(\frac{1}2c-a_h\right)v_j-a_jz-\sum\limits_{ 1\leq k<j}d_{{k,j}}v_k+\sum\limits_{ j<l\leq n}d_{j,l}v_l.\\
\end{split}
\end{equation}

On the other hand, based on (\ref{asos}), we calculate the following equality:
\begin{equation}\begin{split}\label{f2}
\Delta(x)&=\Delta(f-\frac{1}{2}z+v_j)=\Delta(f)-\Delta\left(\frac{1}{2}z\right)+\Delta(v_j)=-2\alpha_{h}^{(f)}f-\sum\limits_{k=1}^n\alpha^{(f)}_{u_k}v_k+\\
&+\left(\frac{\lambda^{(v_j)}}{2}-\alpha^{(v_j)}_{h}\right)v_j-\sum\limits_{ 1\leq k<j}\mu^{(v_j)}_{{k,j}}v_k+\sum\limits_{ j<l\leq n}\mu^{(v_j)}_{{j,l}}v_l.\\
\end{split}
\end{equation}
    Comparing the coefficients at the basis elements $f,\ z$ and $v_j$, (\ref{f1}) and (\ref{f2}) we abtain that
\begin{equation}\label{f3}
    \alpha_h^{(f)}=\alpha_h^{(v_j)}+\alpha_{u_j}^{(f)}-\frac{\lambda^{(v_j)}}{2}.
\end{equation}

 From the equalities, take an element $y=f+\frac{1}{2}z-v_j$ for any fixed $1\leq j\leq n,$

\begin{equation}\begin{split}\label{f11}
\Delta(y)&=\Delta(f+\frac{1}{2}z-v_j)=[b,f+\frac{1}{2}z-v_j]+c\sigma(f+\frac{1}{2}z-v_j)+\\
&+\sum\limits_{1\leq k< l\leq n}d_{k,l}\sigma_{k,l}(f+\frac{1}{2}z-v_j)=
             a_{e}h-2a_{h}f-\sum\limits_{k=1}^{n}a_{k}v_{k}+\frac{c}{2}z-\\
             &-a_eu_j-\left(\frac{1}2c-a_h\right)v_j+a_jz+\sum\limits_{ 1\leq k<j}d_{{k,j}}v_k-\sum\limits_{ j<l\leq n}d_{j,l}v_l.\\
\end{split}
\end{equation}

On the other hand, based on (\ref{asos}), we calculate the following equality:
\begin{equation}\begin{split}\label{f22}
\Delta(y)&=\Delta(f+\frac{1}{2}z-v_j)=\Delta(f)+\Delta\left(\frac{1}{2}z\right)-\Delta(v_j)=-2\alpha_{h}^{(f)}f-\sum\limits_{k=1}^n\alpha^{(f)}_{u_k}v_k-\\
&-\left(\frac{\lambda^{(v_j)}}{2}-\alpha^{(v_j)}_{h}\right)v_j+\sum\limits_{ 1\leq k<j}\mu^{(v_j)}_{{k,j}}v_k-\sum\limits_{ j<l\leq n}\mu^{(v_j)}_{{j,l}}v_l.\\
\end{split}
\end{equation}
    Comparing the coefficients at the basis elements $f,\ z$ and $v_j$, (\ref{f11}) and (\ref{f22}) we obtain that
\begin{equation}\label{f33}
    \alpha_h^{(f)}=\alpha_h^{(v_j)}-\alpha_{u_j}^{(f)}-\frac{\lambda^{(v_j)}}{2}.
\end{equation}

Comparing (\ref{f3}) and (\ref{f33}) we obtain that
$$\alpha_{u_j}^{(f)}=0,\quad \alpha_h^{(f)}=\alpha_h^{(v_j)}-\frac{\lambda^{(v_j)}}{2}.$$

From the equalities, take an element $x=e+\frac{1}{2}z+u_j$, for any fixed $1\leq j\leq n,$

\begin{equation}\begin{split}\label{e1}
\Delta(x)&=\Delta(e+\frac{1}{2}z+u_j)=[b,e+\frac{1}{2}z+u_j]+c\sigma(e+\frac{1}{2}z+u_j)+\\
&+\sum\limits_{1\leq k< l\leq n}d_{k,l}\sigma_{k,l}(e+\frac{1}{2}z+u_j)=
             a_{f}h+2a_he-\sum\limits_{k=1}^{n}b_{k}u_k+\frac{c}{2}z+\\
             &+a_fv_j+\left(\a_h+\frac{1}2c\right)u_j-b_jz-\sum\limits_{ 1\leq k<j}d_{k,j}u_k+\sum\limits_{ j<l\leq n}d_{{j,l}}u_l.\\
\end{split}
\end{equation}

On the other hand, based on (\ref{asos}), we calculate the following equality:
\begin{equation}\begin{split}\label{e2}
\Delta(x)&=\Delta(e+\frac{1}{2}z+u_j)=\Delta(e)+\Delta\left(\frac{1}{2}z\right)+\Delta(u_j)=2\alpha_h^{(e)}e-\sum\limits_{k=1}^n\alpha^{(e)}_{v_k}u_k+\\
&+\left(\alpha^{(u_j)}_{h}+\frac{\lambda^{(u_j)}}{2}\right)u_j-\sum\limits_{ 1\leq k<j}\mu^{(u_j)}_{{k,j}}u_k+\sum\limits_{ j<l\leq n}\mu^{(u_j)}_{{j,l}}u_l.
\end{split}
\end{equation}
    Comparing the coefficients at the basis elements $e,\ z$ and $u_j$, (\ref{e1}) and (\ref{e2}) we obtain that
\begin{equation}\label{e3}
    \alpha_h^{(e)}=\alpha_h^{(u_j)}-\alpha_{v_j}^{(e)}+\frac{\lambda^{(u_j)}}{2}.
\end{equation}

From the relations, take an element $y=e+\frac{1}{2}z-u_j$, for any fixed $1\leq j\leq n,$

\begin{equation}\begin{split}\label{e11}
\Delta(y)&=\Delta(e+\frac{1}{2}z-u_j)=[b,e+\frac{1}{2}z-u_j]+c\sigma(e+\frac{1}{2}z-u_j)+\\
&+\sum\limits_{1\leq k< l\leq n}d_{k,l}\sigma_{k,l}(e+\frac{1}{2}z-u_j)=
             a_{f}h+2a_he-\sum\limits_{k=1}^{n}b_{k}u_k+\frac{c}{2}z-\\
             &-a_fv_j-\left(\a_h+\frac{1}2c\right)u_j+b_jz+\sum\limits_{ 1\leq k<j}d_{k,j}u_k-\sum\limits_{ j<l\leq n}d_{{j,l}}u_l.\\
\end{split}
\end{equation}

On the other hand, based on (\ref{asos}), we calculate the following equality:
\begin{equation}\begin{split}\label{e22}
\Delta(y)&=\Delta(e+\frac{1}{2}z-u_j)=\Delta(e)+\Delta\left(\frac{1}{2}z\right)-\Delta(u_j)=2\alpha_h^{(e)}e-\sum\limits_{k=1}^n\alpha^{(e)}_{v_k}u_k-\\
&-\left(\alpha^{(u_j)}_{h}+\frac{\lambda^{(u_j)}}{2}\right)u_j+\sum\limits_{ 1\leq k<j}\mu^{(u_j)}_{{k,j}}u_k-\sum\limits_{ j<l\leq n}\mu^{(u_j)}_{{j,l}}u_l.
\end{split}
\end{equation}
    Comparing the coefficients at the basis elements $e,\ z$ and $u_j$, (\ref{e11}) and (\ref{e22}) we obtain that
\begin{equation}\label{e33}
    \alpha_h^{(e)}=\alpha_h^{(u_j)}+\alpha_{v_j}^{(e)}+\frac{\lambda^{(u_j)}}{2}.
\end{equation}

Comparing (\ref{e3}) and (\ref{e33}) we get that
$$\alpha_{v_j}^{(e)}=0,\quad \alpha_h^{(e)}=\alpha_h^{(u_j)}+\frac{\lambda^{(u_j)}}{2}.$$

\end{proof}

\begin{lemma}\label{lem6}
\begin{equation}
\begin{aligned}
\Delta(u_j)&=\alpha_h^{(e)}u_j+\sum\limits_{ 1\leq k<j}\mu_{{k,j}}u_k-\sum\limits_{ j<l\leq n}\mu_{{j,l}}u_l,\\
\Delta(v_j)&=-\alpha_h^{(e)}v_j+\sum\limits_{ 1\leq k<j}\mu_{{k,j}}v_k-\sum\limits_{ j<l\leq n}\mu_{{j,l}}v_l.\\
\end{aligned}
\end{equation}
\end{lemma}

\begin{proof}

From the equalities, take an element $x=u_p+iu_j$, for any fixed $1\leq p<j\leq n,$
and $i$ complex number
\begin{equation}\begin{split}\label{u1}
\Delta(x)&=\Delta(u_p+iu_j)=[b,u_p+iu_j]+c\sigma(u_p+iu_j)+
\sum\limits_{1\leq k< l\leq n}d_{k,l}\sigma_{k,l}(u_p+iu_j)=\\
&=a_fv_p+\left(\a_h+\frac{1}2c\right)u_p-b_pz-\sum\limits_{ 1\leq k<p}d_{k,p}u_k+\sum\limits_{ p<l\leq n}d_{{p,l}}u_l+\\
&+ia_fv_j+i\left(\a_h+\frac{1}2c\right)u_j-ib_jz-i\sum\limits_{ 1\leq k<j}d_{k,j}u_k+i\sum\limits_{ j<l\leq n}d_{{j,l}}u_l.\\
\end{split}
\end{equation}

On the other hand, based on (\ref{asos}), we calculate the following equality:
\begin{equation}\begin{split}\label{u2}
\Delta(x)&=\Delta(u_p+iu_j)=\Delta(u_p)+i\Delta(u_j)=\alpha_h^{(e)}u_p-\sum\limits_{ 1\leq k<p}\mu^{(u_p)}_{{k,p}}u_k-\\
&+\sum\limits_{ p<l\leq n}\mu^{(u_p)}_{{p,l}}u_l+
i\alpha_h^{(e)}u_j-i\sum\limits_{ 1\leq k<j}\mu^{(u_j)}_{{k,j}}u_k+i\sum\limits_{ j<l\leq n}\mu^{(u_j)}_{{j,l}}u_l.\\
\end{split}
\end{equation}
     Comparing the coefficients at the basis elements $u_p$ and $u_j$, (\ref{u1}) and (\ref{u2}) we get
\begin{equation*}
    (a_h+\frac{1}2c)-id_{p,j}=\alpha_h^{(e)}-i\mu^{(u_j)}_{{p,j}},\\ \quad
    i(a_h+\frac{1}2c)+d_{p,j}=i\alpha_h^{(e)}+\mu^{(u_p)}_{{p,j}},
\end{equation*}
which implies
$$\mu^{(u_p)}_{{p,j}}=\mu^{(u_j)}_{{p,j}}.$$

Now, if we look at the above element $y=v_p+iv_j,$ for any fixed $1\leq p<j\leq n,$, the following result is obtained
$$\mu^{(v_p)}_{{p,j}}=\mu^{(v_j)}_{{p,j}}.$$

\end{proof}

\begin{lemma}\label{lem81}

$\Delta$ is a derivation on $\mathcal{S}.$

%\begin{equation}
%\begin{aligned}
%\Delta(e)&=2\alpha_he,\\
%\Delta(h)&=0,\\
%\Delta(f)&=-2\alpha_{h}f,\\
%\Delta(u_j)&=\alpha_{h}u_j-\sum\limits_{ 1\leq k<j}\mu_{{k,j}}u_k+\sum\limits_{ j<l\leq n}\mu_{{j,l}}u_l,\\
%\Delta(v_j)&=\alpha_{h}v_j-\sum\limits_{ 1\leq k<j}\mu_{{k,j}}v_k+\sum\limits_{ j<l\leq n}\mu_{{j,l}}v_l,\\
%\Delta(z)&=0.
%\end{aligned}
%\end{equation}
\end{lemma}

\begin{proof}
    For any element $x=x_ff+x_hh+x_ee+x_zz+\sum\limits_{j=1}^nx_{u_j}u_j+\sum\limits_{j=1}^nx_{v_j}v_j,$ such that
\begin{equation*}\begin{split}
\Delta(x)&=-2\alpha_hx_ff+2\alpha_hx_ee+
\sum\limits_{j=1}^n\big(\alpha_{h}x_{u_j}u_j-\sum\limits_{1\leq k<j}\mu_{{k,j}}x_{u_k}u_k+\sum\limits_{ j<l\leq n}\mu_{{j,l}}x_{u_l}u_l+\\
&+\alpha_{h}x_{v_j}v_j-\sum\limits_{1\leq k<j}\mu_{{k,j}}x_{v_k}v_k+\sum\limits_{j<l\leq n}\mu_{{j,l}}x_{v_l}v_l\big)=
[a_hh,x]+ \sum\limits_{1\leq k<l\leq n}\mu_{k,l}\sigma_{k,l}(x).
\end{split}\end{equation*}

Then, $\Delta$ is a derivation.
\end{proof}

Now we are in position to prove Theorem \ref{thm21}.

\begin{proof}[Proof of Theorem \ref{thm21}] Let $\nabla$ be a local derivation of $\mathcal{S}_n.$ For the element $h+z$ exists a derivation
$D_{h+z}$ such that
$$\nabla(h+z)=D_{h+z}(h+z).$$

Set $\Delta=\nabla-D_{h+z}.$  Then $\Delta$ is a local derivation such that $\Delta(h+z)=0.$
By Lemma \ref{lem81}, $\Delta(x) = [a_hh,x]+ \sum\limits_{1\leq k<l\leq n}\mu_{k,l}\sigma_{k,l}(x).$
Thus
$$\nabla(x)=[a_hh,x]+ \sum\limits_{1\leq k<l\leq n}\mu_{k,l}\sigma_{k,l}(x)+D_{h+z}(x).$$  Then $\nabla$ is a derivation. The proof is complete.
\end{proof}

\begin{rem}
    It should be noted that, if $n = 1$, in \cite{AlaYus2} it is proved every local derivation is a derivation.
\end{rem}

\textbf{Data availability}  There is no other data for this article

\end{document}